\documentclass[11pt]{amsart}
\usepackage{amssymb, latexsym}
\theoremstyle{plain}
\newtheorem{theorem}{Theorem}
\newtheorem{corollary}{Corollary}

\newtheorem*{2'}{Theorem 2'}
\newtheorem*{3'}{Theorem 3'}

\theoremstyle{remark}

\newtheorem*{Remark 1}{Remark 1}
\newtheorem*{Remark 2}{Remark 2}
\newtheorem*{Remark 3}{Remark 3}
\newtheorem*{Remark 4}{Remark 4}

\numberwithin{equation}{section}

\begin{document}

\title [Probabilistic Proofs of Some Generalized Mertens' Formulas ]
{Probabilistic Proofs of Some Generalized Mertens' Formulas Via Generalized Dickman Distributions }

\author{Ross G. Pinsky}


\address{Department of Mathematics\\
Technion---Israel Institute of Technology\\
Haifa, 32000\\ Israel}
\email{ pinsky@math.technion.ac.il}

\urladdr{http://www.math.technion.ac.il/~pinsky/}

\date{}

\begin{abstract}

The classical Mertens' formula states that
$
\prod_{p\le N}\big(1-\frac1p)^{-1}\sim e^\gamma\log N,
$
where the product is over all primes $p$ less than or equal to $N$, and $\gamma$ is the Euler-Mascheroni constant.
By the Euler product formula,  this is equivalent to either of the following statements:
$$
\begin{aligned}
&i. \lim_{N\to\infty}\frac{\sum_{n:p|n\Rightarrow p\le N}\thinspace\frac1n}{\sum_{n\le N}\frac1n}=e^\gamma\ \
&ii. \sum_{n:p|n\Rightarrow p\le N}\thinspace\frac1n\sim e^\gamma\log N.
\end{aligned}
$$
Via some random integer constructions and a criterion for weak convergence of distributions to so-called generalized Dickman distributions,
we obtain some generalized Mertens' formulas, some of which are new and some of which have been proved using number-theoretic tools.
For example, in the spirit of (i), we show that if $A$ is a subset of the primes which has
natural density $\theta\in(0,1]$ with respect to the set of all primes, then
$$
\lim_{N\to\infty}\frac{\sum_{n:p|n\Rightarrow p\le N\thinspace\text{and}\thinspace p\in A}\frac1n}
{\sum_{n\le N:p|n\Rightarrow p\in A}\frac1n}=e^{\gamma\theta}\Gamma(\theta+1),
$$
and also,  for any $k\ge2$,
$$
\lim_{N\to\infty}\frac{\sum^{'(k)}_{n:p|n\Rightarrow p\le N\thinspace\text{and}\thinspace p\in A}\frac1n}
{\sum^{'(k)}_{n\le N:p|n\Rightarrow p\in A}\frac1n}=e^{\gamma\theta}\Gamma(\theta+1),
$$
where $\sum^{'(k)}$ denotes  that the summation is restricted to  $k$-free positive integers.
In the spirit of (ii), we show for example that

$
\sum^{'(k)}_{n:p|n\Rightarrow p\le N}\frac1{n_{\{(k-1)-\text{free}\}}\phi(n_{\{(k-1)-\text{power}\}})}\sim e^\gamma\log N,
$
\medskip

\noindent where $\phi$ is the Euler totient function, and $n_{\{(k-1)-\text{free}\}}$ and  $n_{\{(k-1)-\text{power}\}}$ are the $(k-1)$-free part and the $(k-1)$-power part of $n$.
    \end{abstract}

\maketitle
\section{Introduction and Statement of Results}\label{intro}
The classical Mertens' formula states that
\begin{equation}\label{Mertens}
\prod_{p\le N}\big(1-\frac1p)^{-1}\sim e^\gamma\log N\approx 1.78\log N,\ \text{as}\ N\to\infty,
\end{equation}
where the product is over all primes $p$ less than or equal to $N$, and $\gamma$ is the Euler-Mascheroni constant.
(Actually, the classical formula states a little more; namely,
$\prod_{p\le N}\big(1-\frac1p)^{-1}= e^\gamma\log N+O(1)$.)
By the Euler product formula,  \eqref{Mertens} is equivalent to either of the following statements, which are more in the spirit
 of the results we present in this paper:
\begin{equation}\label{Mertens-alt}
\begin{aligned}
&i. \lim_{N\to\infty}\frac{\sum_{n\in\mathbb{N}:p|n\Rightarrow p\le N}\thinspace\frac1n}{\sum_{n\le N}\frac1n}=e^\gamma;\\
&ii. \sum_{n\in\mathbb{N}:p|n\Rightarrow p\le N}\thinspace\frac1n\sim e^\gamma\log N,\ \text{as}\ N\to\infty,
\end{aligned}
\end{equation}
where as usual, $[N]$ denotes the set $\{1,\cdots, N\}$.
In words, (\ref{Mertens-alt}-i) states that $e^\gamma$ is the limit of  the ratio between the harmonic series restricted to the positive integers all of whose prime factors are no greater than $N$ and the harmonic series
restricted to the positive integers no greater than $N$.

One of the results in  \cite{P19} involved the construction of a sequence of random integers whose    distributions were shown to converge weakly to the so-called Dickman distribution. It was noted
in that paper that Mertens' formula   follows readily as a corollary of this result.
In this paper, we make a number of random integer  constructions in a similar vein,  and use our recent paper \cite{P18} to show that their distributions converge weakly to so-called generalized Dickman distributions.
From these results, we obtain  a number of generalizations of Mertens' formula,
some of which are known via number theoretic methods, and some of which appear to be new.

We begin by
 introducing  some notation and constructing the three sequences of random integers  that will be used in this paper.
 Then we state our generalized Mertens' formulas.

Let $\mathbb{P}$ denote the set of prime numbers.
Recall that for $k\ge2$,  an integer $n\in\mathbb{N}$ is called \it\ $k$-free\rm\ if $p^k\nmid n$, for all primes $p$.
Let $A\subset\mathbb{P}$ be an infinite set of primes.
Denote the primes in $A$ in increasing order by $p_{1;A},p_{2;A},\cdots$.
Let $\{T_j\}_{j=1}^\infty$ be a sequence of independent random variables with $T_j$ distributed according to the geometric distribution with parameter $\frac1{p_{j;A}}$
$\big(T_j\sim\text{Geom}(\frac1{p_{j;A}})\big)$, $j=1,\cdots$;
that is
\begin{equation}\label{Tdist}
P(T_j=m)=(1-\frac1{p_{j;A}})(\frac1{p_{j;A}})^m,\ m=0,1,\cdots.
\end{equation}
For  $N\in\mathbb{N}$, we define a random integer by
\begin{equation}\label{randomint1}
I_{N;A,1}=\prod_{j=1}^Np_{j;A}^{T_j}.
\end{equation}
By construction, the support of  $I_{N;A,1}$ is $\{n\in\mathbb{N}: p|n\Rightarrow p\le N\thinspace\text{and}\thinspace p\in A\}$.
See \cite{P19} for a detailed study of the random integer sequence $\{I_{N;A,1}\}_{N=1}^\infty$ when $A=\mathbb{P}$.

Let $k\ge2$.
We define a second random integer sequence by replacing the random variables $\{T_j\}_{j=1}^\infty$ by
a sequence $\{U_j\}_{j=1}^\infty$ of independent random variables, where
$U_j$ is distributed as $T_j$ conditioned on being less than $k$
$\big(U_j\stackrel{\text{dist}}{=}T_j|\{T_j<k\}\big)$;  that is
\begin{equation}\label{Udist}
P(U_j=m)=P(T_j=m|T_j<k)=\frac{1-\frac1{p_{j;A}}}{1-(\frac1{p_{j;A}})^k}(\frac1{p_{j;A}})^m,\ m=0,1,\cdots, k-1.
\end{equation}
For  $N\in\mathbb{N}$, define a random integer by
\begin{equation}\label{randomint2}
I_{N;A,2}=\prod_{j=1}^Np_{j;A}^{U_j}.
\end{equation}
By construction, the support of  $I_{N;A,2}$ is the set of $k$-free integers in $\{n\in\mathbb{N}: p|n\Rightarrow p\le N\thinspace\text{and}\thinspace p\in A\}$.

Finally, for $k\ge2$, we construct a third random integer sequence from a sequence $\{V_j\}_{j=1}^\infty$
of independent random variables, where
$V_j$ is distributed as $T_j$ truncated at $k-1$
$\big(V_j\stackrel{\text{dist}}{=}T_j\wedge (k-1)\big)$;  that is
\begin{equation}\label{Vdist}
\begin{aligned}
&P(V_j=m)=(1-\frac1{p_{j;A}})(\frac1{p_{j;A}})^m,\ m=0,\cdots, k-2;\\
&P(V_j=k-1)=(\frac1{p_{j;A}})^{k-1}.
\end{aligned}
\end{equation}
For  $N\in\mathbb{N}$, define a random integer by
\begin{equation}\label{randomint3}
I_{N;A,3}=\prod_{j=1}^Np_{j;A}^{V_j}.
\end{equation}
By construction, the support of $I_{N;A,3}$ is the set of $k$-free integers in $\{n\in\mathbb{N}: p|n\Rightarrow p\le N\thinspace\text{and}\thinspace p\in A\}$.

We now present our generalized Mertens' formulas.
For  $A\subset\mathbb{P}$, denote by
$$
D_{\text{nat-prime}}(A):=\lim_{N\to\infty}\frac{|A\cap[N]|}{|\mathbb{P}\cap[N]|}
$$
the natural density of $A$ in $\mathbb{P}$, if it exists.
For $B\subset\mathbb{N}$, let $\sum^{'(k)}_B$ denote the summation restricted to the $k$-free powers in $B$.
Let $\Gamma(x):=\int_0^\infty t^{x-1}e^{-t}dt,\ x>0$, denote the Gamma function.
\begin{theorem}\label{GenMert}
Let $A\subset\mathbb{P}$ be a subset of primes whose natural density in $\mathbb{P}$ is $D_{\text{nat-prime}}(A)=\theta\in(0,1]$.  Then

\noindent i.
\begin{equation}\label{basic}
\lim_{N\to\infty}\frac{\sum_{n\in\mathbb{N}:p|n\Rightarrow p\le N\thinspace\text{and}\thinspace p\in A}\thinspace\frac1n}{\sum_{n\le N:p|n\Rightarrow p\in A}\frac1n}=e^{\gamma\theta}\thinspace \Gamma(\theta+1);
\end{equation}

\noindent ii. for  $k\ge2$,
$$
\lim_{N\to\infty}\frac{\sum^{'(k)}_{n\in\mathbb{N}:p|n\Rightarrow p\le N\thinspace\text{and}\thinspace p\in A}\thinspace\frac1n}
{\sum^{'(k)}_{n\le N:p|n\Rightarrow p\in A}\frac1n}=e^{\gamma\theta}\thinspace\Gamma(\theta+1).
$$
\end{theorem}

\noindent \bf Remark 1.\rm\ Part (i) of Theorem \ref{GenMert} in  a slightly different but equivalent
form appears in \cite{Wirsing}, and a refined version appears in \cite{R} .
Part (ii) seems to be new.

\noindent\bf Remark 2.\rm\ The function $\theta\to e^{\gamma\theta}\Gamma(\theta+1)$, $\theta\in(0,1]$, is increasing, is equal to 1 at $\theta=0^+$ and is equal to $e^{\gamma}$ at $\theta=1$.

\noindent \bf Remark 3.\rm\ When $A=\mathbb{P}$, (i) reduces to the classical Mertens' formula.

\noindent \bf Remark 4.\rm\ For $l\in\mathbb{N}$ and $j$ satisfying $1\le j<l$ and $(j,l)=1$, let $A_{l;j}=\{p\in\mathbb{P}: p=j\thinspace \text{mod}\thinspace l\}$ denote the set of primes that are equal
to $j$ modulo $l$. Dirichlet's arithmetic progression theorem states that  $D_{\text{nat-prime}}(A_{l;j})=\frac1{\phi(l)}$, where $\phi$ is Euler's totient function.
In \cite{W}, it was proved that
\begin{equation}\label{Williams}
\prod_{p\le N,p=j\thinspace \text{mod}\thinspace l}(1-\frac1p)^{-1}=
\sum_{n:p|n\Rightarrow p\le N\thinspace\text{and}\thinspace p\in A_{l;j}}\thinspace\frac1n\sim \Big(e^{\gamma}\frac{\phi(l)}{l}C(l,j)\Big)^\frac1{\phi(l)}(\log N)^{\frac1{\phi(l)}},
\end{equation}
as $N\to\infty$, where $C(l,j)$ is a complicated expression involving Dirichlet characters modulo $l$. Thus, by part (i) we obtain
$$
\sum_{n\le N: p|n\Rightarrow p\in A_{l;j}}\frac1n\sim\frac1{\Gamma(1+\frac1{\phi(l)})}\Big(\frac{\phi(l)}{l}C(l,j)\Big)^\frac1{\phi(l)}(\log N)^{\frac1{\phi(l)}},\ \text{as}\ N\to\infty,
$$
as was noted in \cite{W}.
A much simpler looking form for $C(l,j)$ was obtained in \cite{LZ1}; namely,
$$
\begin{aligned}
&\frac {\phi(l)}lC(l,j)=\prod_p(1-\frac1p)^{-\alpha(p,l,j)},\ \text{where}\ \alpha(p,l,j)=\begin{cases}\phi(l)-1,\ p=j\ \text{mod}\ l;\\ -1,\ \text{otherwise}.\end{cases}
\end{aligned}
$$
See also \cite{LZ2}  for more on this constant.

\noindent From part (ii), we obtain for $k\ge2$,
$$
\lim_{N\to\infty}\frac{\sum^{'(k)}_{n:p|n\Rightarrow p\le N\thinspace\text{and}\thinspace p\in A_{l;j}}\thinspace\frac1n}
{\sum^{('(k)}_{n\le N:p|n\Rightarrow p\in A_{l;j}}\frac1n}=e^{\frac{\gamma}{\phi(l)}}\thinspace\Gamma(\frac1{\phi(l)}+1).
$$

\bigskip

Whereas Theorem \ref{GenMert} is an asymptotic result and depends on the above-mentioned convergence to generalized
Dickman distributions, the next result is a non-asymptotic identity that holds for all $N$, and that only requires
the random integer  $I_{N;A,3}$.
For each $k\ge2$, every  $n\in\mathbb{N}$ can be written uniquely as $n=n_{\{k-\text{free\}}}n_{\{k-\text{power}\}}$, where
$n_{\{k-\text{free\}}}$ is $k$-free and   $n_{\{k-\text{power}\}}$ is a $k$th power. We  extend this to $k=1$ by defining $n_{\{1-\text{free}\}}=1$ and
$n_{\{1-\text{power}\}}=n$.
\begin{theorem}\label{fixedN}
Let $N\ge1$ and let $A_N=\{p_{1;A},\cdots, p_{N;A}\}\subset \mathbb{P}$. Then for $k\ge2$,

$
\sum^{'(k)}_{n:p|n\Rightarrow p\in A_N}\thinspace\frac1{n_{\{(k-1)-\text{free}\}}\phi(n_{\{(k-1)-\text{power}\}})}=
\sum_{n:p|n\Rightarrow p\in A_N}\frac1n,
$
\medskip

\noindent where $\phi$ is the Euler totient fuction.
\end{theorem}
\bf\noindent Remark.\rm\ Theorem \ref{fixedN} seems to be new. Note that when $k=2$, the result is

$
\sum^{'(2)}_{n:p|n\Rightarrow p\in A_N}\thinspace\frac1{\phi(n)}=
\sum_{n:p|n\Rightarrow p\in A_N}\frac1n.
$
\medskip

We have the following corollary.
\begin{corollary}\label{fixedNcor}
For $k\ge2$,
\medskip

$
\sum^{'(k)}_{n:p|n\Rightarrow p\le N}\frac1{n_{\{(k-1)-\text{free}\}}\phi(n_{\{(k-1)-\text{power}\}})}\sim e^\gamma\log N,\ \text{as}\ N\to\infty.
$
\end{corollary}
\noindent \it Proof of Corollary.\rm\
 By Theorem \ref{fixedN} with $A_N$ replaced by $\mathbb{P}\cap[N]$, we have
 \medskip

 $
 \sum^{'(k)}_{n:p|n\Rightarrow p\le N}\thinspace\frac1{n_{\{(k-1)-\text{free}\}}\phi(n_{\{(k-1)-\text{power}\}})}=
\sum_{n:p|n\Rightarrow p\le N}\frac1n,
 $
\medskip

\noindent and by Mertens' formula (\ref{Mertens-alt}-ii) we have
$\sum_{n:p|n\Rightarrow p\le N}\frac1n\sim e^\gamma\log N$.
\hfill$\square$
\medskip

Our final theorem combines some of the ingredients of Theorems \ref{GenMert} and \ref{fixedN}.
\begin{theorem}\label{GenPhiMert}
Let $A\subset\mathbb{P}$ be a subset of primes whose density in $\mathbb{P}$ is $D_{\text{nat-prime}}(A)=\theta\in(0,1]$.  Then
 for  all $k\ge 2$,
\begin{equation}\label{genphimert}
\lim_{N\to\infty}\frac{\sum_{n:p|n\Rightarrow p\le N\thinspace\text{and}\thinspace p\in A}\frac1n}
{\sum^{'(k)}_{n\le N:p|n\Rightarrow p\in A} \thinspace\frac1{n_{\{(k-1)-\text{free}\}}\phi(n_{\{(k-1)-\text{power}\}})}}
=e^{\gamma\theta}\Gamma(\theta+1).
\end{equation}
\end{theorem}
\bf\noindent Remark.\rm\ Theorem \ref{GenPhiMert} seems to be new.
Note that when $k=2$, the result is
$$
\lim_{N\to\infty}\frac{\sum_{n:p|n\Rightarrow p\le N\thinspace\text{and}\thinspace p\in A}\thinspace\frac1n}
{\sum^{'(2)}_{n\le N:p|n\Rightarrow p\in A}\frac1{\phi(n)}}=e^{\gamma\theta}\thinspace\Gamma(\theta+1).
$$

\medskip

From Theorem \ref{GenPhiMert}, we obtain the following corollary
\begin{corollary}\label{totientharmcor}
Let $A\subset\mathbb{P}$ be a subset of primes whose density in $\mathbb{P}$ is $D_{\text{nat-prime}}(A)=\theta\in(0,1]$.  Then
 for  all $k\ge 2$,
\medskip

$
\sum^{'(k)}_{n\le N:p|n\Rightarrow p\in A} \frac1{n_{\{(k-1)-\text{free}\}}\phi(n_{\{(k-1)-\text{power}\}})}\sim
\sum_{n\le N: p|n\Rightarrow p\in A}\frac1n,\ \text{as}\ N\to\infty.
$
\end{corollary}
\noindent \it Proof of Corollary.\rm\ Compare \eqref{basic} to \eqref{genphimert}.
\hfill $\square$

\noindent \bf Remark.\rm\ When $A=\mathbb{P}$, Corollary \ref{totientharmcor} reduces to

$
\sum^{'(k)}_{n\le N:}\frac1{n_{\{(k-1)-\text{free}\}}\phi(n_{\{(k-1)-\text{power}\}})}\sim\log N,\ \text{as}\ N\to\infty, \ \text{for all}\ k\ge2.
$
\medskip

\noindent When $k=2$, this reduces to
\medskip

$
\sum^{'(2)}_{n\le N}\frac1{\phi(n)}\sim\log N,\ \text{as}\ N\to\infty,
$
\medskip

\noindent which is known (see \cite{Ward} or \cite[p. 43, problem 17]{MV}.
The generalization  to all $k\ge2$ seems to be new.
As an aside, we note that
$$
\sum_{n\le N}\frac1{\phi(n)}\sim\frac{\zeta(2)\zeta(3)}{\zeta(6)}\log N\approx 1.94\log N,\ \text{as}\ N\to\infty.
$$
(see \cite[p. 42, problem 13-(d)]{MV}).


We prove Theorems \ref{GenMert}--\ref{GenPhiMert}
in sections \ref{2}--\ref{4} respectively.

\section{Proof of Theorem \ref{GenMert}}\label{2}
Fix a subset $A\subset\mathbb{P}$ which satisfies $D_{\text{nat-prime}}(A)=\theta\in(0,1]$.
Denote the primes in $A$ in increasing order by $p_{1;A},p_{2;A},\cdots$, and let
$$
A_N=\{p_{1;A},p_{2;A},\cdots, p_{N;A}\}.
$$
\noindent \it Proof of part (i).\rm\
Let the random integer $I_{N;A,1}$ be as in \eqref{randomint1}, where $\{T_j\}_{j=1}^\infty$ is a sequence of independent random variables with distributions given by \eqref{Tdist}.
The support of  $I_{N;A,1}$ is $\{n\in\mathbb{N}:p|n\Rightarrow p\in A_N\}$, and for  arbitrary $n=\prod_{j=1}^Np_{j;A}^{c_j}$
in the support,
\begin{equation}\label{1nprob}
\begin{aligned}
&P(I_{N;A,1}=n)=P(T_j=c_j,\thinspace j\in[N])=\prod_{j=1}^NP(T_j=c_j)\\
&=\prod_{j=1}^N(1-\frac1{p_{j;A}})(\frac1{p_{j;A}})^{c_j}=\frac1n\prod_{j=1}^N(1-\frac1{p_{j;A}}).
\end{aligned}
\end{equation}

We have
\begin{equation*}
\log I_{N;A,1}=\sum_{j=1}^NT_j\log p_{j;A}.
\end{equation*}
Noting that the expected value of $T_j$ is given by
\begin{equation}\label{ET}
ET_j=\frac1{p_{j;A}-1},
\end{equation}
we have
$$
E\log I_{N;A,1}=\sum_{j=1}^N\frac{\log p_{j;A}}{p_{j;A}-1}.
$$
It follows by the assumption on the density of $A$ and  by the prime number theorem that
\begin{equation}\label{pjAasymp}
p_{j;A}\sim\frac {j\log j}\theta,\ \text{as}\ j\to\infty,
\end{equation}
and thus that
\begin{equation}\label{Elog}
E\log I_{N;A,1}\sim\theta\log N,\ \text{as}\ N\to\infty.
\end{equation}

We will demonstrate below that the conditions of a theorem in \cite{P18} are satisfied, from which it follows
that
\begin{equation}\label{strangethm}
\lim_{N\to\infty}\frac{\log I_{N;A,1}}{E\log I_{N;A,1}}\stackrel{\text{dist}}{=}\frac1\theta D_\theta,
\end{equation}
where $D_\theta$ is a random variable distributed according to  the generalized Dickman distribution GD$(\theta)$ with
parameter $\theta$. This distribution has density function
$p_\theta=\frac{e^{-\gamma\theta}}{\Gamma(\theta)}\rho_\theta$, where  $\rho_\theta$ satisfies the differential-delay equation
\begin{equation}\label{rhotheta}
\begin{aligned}
&\rho_\theta(x)=0,\ x\le0;\\
&\rho_\theta(x)=x^{\theta-1},\ 0< x\le 1;\\
&x\rho_\theta'(x)+(1-\theta)\rho_\theta(x)+\theta\rho_\theta(x-1)=0,\ x>1.
\end{aligned}
\end{equation}
(The function $\rho_1$ is known as the Dickman function; we call $\rho_\theta$ a generalized Dickman function.)

On the one hand, by the convergence in distribution in \eqref{strangethm} and the fact that the limiting distribution is a continuous one, for any
sequence $\{\theta_N\}_{N=1}^\infty$ satisfying $\lim_{N\to\infty}\theta_N=\theta$, we have
\begin{equation}\label{onehand}
\begin{aligned}
&\lim_{N\to\infty}P(\frac{\log I_{N;A,1}}{E\log I_{N;A,1}}\le \frac1{\theta_N})=P(\frac1\theta D_\theta\le \frac1\theta)
=\int_0^1p_\theta(x)dx\\
&=\frac{e^{-\gamma\theta}}{\Gamma(\theta)}\int_0^1x^{\theta-1}dx=
\frac{e^{-\gamma\theta}}{\theta\Gamma(\theta)}=
\frac{e^{-\gamma\theta}}{\Gamma(\theta+1)}.
\end{aligned}
\end{equation}
On the other hand,
let $\theta_N:=\frac{E\log I_{N;A,1}}{\log p_{N;A}}$ and note from \eqref{pjAasymp} and \eqref{Elog} that
$\lim_{N\to\infty}\theta_N=\theta$.
It follows from
 \eqref{1nprob} that
\begin{equation}\label{otherhand}
\begin{aligned}
&P(\frac{\log I_{N;A,1}}{E\log I_{N;A,1}}\le \frac1{\theta_N})=P\big(I_{N;A,1}\le \exp(\frac{E\log I_{N;A,1}}{\theta_N})\big)\\
&=P(I_{N;A,1}\le p_{N;A})=
\Big(\prod_{j=1}^N(1-\frac1{p_{j;A}})\Big)\sum_{n\le p_{N;A}:p|n\Rightarrow p\in A  }\frac1n\\
&=\frac{\sum_{n\le p_{N;A}:p|n\Rightarrow p\in A}\frac1n}{\sum_{n: p|n\Rightarrow p\in A_N}\frac1n}.
\end{aligned}
\end{equation}
From \eqref{onehand} and \eqref{otherhand}, we conclude that
\begin{equation}\label{finalalmost}
\lim_{N\to\infty}\frac{\sum_{n:p|n\Rightarrow p\in A_N}\frac1n}
{\sum_{n\le p_{N;A}:p|n\Rightarrow p\in A}\frac1n}=e^{\gamma\theta}\Gamma(\theta+1).
\end{equation}

Now \eqref{finalalmost}
 is equivalent to part (i) of Theorem \ref{GenMert}.
Indeed, for any $M\in\mathbb{N}$, let $N^+(M)=\max\{n:p_{n;A}\le M\}$.
Then
\begin{equation}\label{equiv1}
\sum_{n:p|n\Rightarrow p\le M\thinspace\text{and}\thinspace p\in A}\frac1n=\sum_{n:p|n\Rightarrow p\le p_{N^+(M);A}\thinspace\text{and}\thinspace p\in A}\frac1n,
\end{equation}
and
\begin{equation}\label{equiv2}
\sum_{n\le M:p|n\Rightarrow p\in A}\frac1n=\sum_{n\le p_{N^+(M);A}:p|n\Rightarrow p\in A}\frac1n+H_M,
\end{equation}
where
\begin{equation}\label{equiv3}
H_M:=\sum_{n\in[p_{N^+(M);A}+1,M]:p|n\Rightarrow p\in A}\frac1n\le
\sum_{n=p_{N^+(M);A}+1}^{p_{N^+(M)+1;A}}\frac1n\le\log\frac{p_{N^+(M)+1;A}}{p_{N^+(M);A}}=o(1).
\end{equation}
Thus, \eqref{finalalmost}-\eqref{equiv3} give
$$
\lim_{M\to\infty}\frac{\sum_{n:p|n\Rightarrow p\le M\thinspace\text{and}\thinspace p\in A}\frac1n}
{\sum_{n\le M:p|n\Rightarrow p\in A}\frac1n}=e^{\gamma\theta}\Gamma(\theta+1),
$$
which is part (i) of the theorem.

We now show that \eqref{strangethm} holds.
Let $\{B_j\}_{j=1}^\infty$ and  $\{Xj\}_{j=1}^\infty$
  be mutually independent random variables with distributions as follows:

\noindent $B_j\sim \text{Ber}(\frac1{p_{j;A}})$; that is,
\begin{equation}\label{qj}
q_j:=P(B_j=1)=1-P(B_j=0)=\frac1p_{j;A}.
\end{equation}

\noindent $X_j\stackrel{\text{dist}}{=}\log p_{j;A}\cdot T_j|\{T_j\ge1\}$;
that is
$$
P(X_j=m\log p_{j;A})=P(T_j=m|T_j\ge1)=
(1-\frac1{p_{j;A}})(\frac1{p_{j;A}})^{m-1},\ m=1,2,\cdots.
$$
Then
\begin{equation}\label{muj}
\mu_j:=EX_j=\frac{p_{j;A}}{p_{j;A}-1}\log p_{j;A},
\end{equation}
and it follows that
\begin{equation*}\label{XEX}
\lim_{j\to\infty}\frac{X_j}{\mu_j}\stackrel{\text{dist}}{=}1.
\end{equation*}
By the construction of $\{B_j\}_{j=1}^\infty$ and  $\{X_j\}_{j=1}^\infty$, we have
\begin{equation}\label{logIBX}
\log I_{N;A,1}=\sum_{j=1}^NT_j\log p_{j;A}\stackrel{\text{dist}}{=}\sum_{j=1}^NB_jX_j.
\end{equation}
From \eqref{pjAasymp}, \eqref{qj} and \eqref{muj},
\begin{equation}\label{qmu}
\mu_j\sim\log j;\ \ \ q_j\sim\frac\theta{j\log j},\ \text{as}\ j\to\infty.
\end{equation}
Let $W_N=\frac{\sum_{j=1}^NB_jX_j}{E\sum_{j=1}^NB_jX_j}$. Now Theorem 1.2 in \cite{P18} applies to $W_N$. Our notation here coincides with  the notation in that theorem
except that the summation there is over $k$ while here it is over $j$, and  $p_k$ there corresponds to $q_j$ here.
In light of \eqref{qmu}, we have $J_\mu=J_p=1,a_0=0,a_1=b_0=b_1=c_\mu=1, c_p=\theta$ in the notation of that theorem.
For these values, the theorem indicates that $W_N$ converges in distribution to $\frac1\theta D_\theta$.
By \eqref{logIBX},  $\frac{\log I_{N;A,1}}{E\log I_{N;A,1}}\stackrel{\text{dist}}{=}W_N$; thus \eqref{strangethm} holds.
\hfill $\square$

\it \noindent Proof of part (ii).\rm\ Fix $k\ge2$. The proof follows the proof of part (i), except that we replace the random integer $I_{N;A,1}$ by the random integer $I_{N;A,2}$ from \eqref{randomint2},
 where $\{U_j\}_{j=1}^\infty$ is a sequence of random variables with distributions given by \eqref{Udist}.
The support of  $I_{N;A,2}$ is the set of $k$-free integers in $\{n\in\mathbb{N}:p|n\Rightarrow p\in A_N\}$,
 and  for arbitrary $n=\prod_{j=1}^Np_{j;A}^{c_j}$ in the support,
\begin{equation}\label{1nprobagain}
\begin{aligned}
&P(I_{N;A,2}=n)=\prod_{j=1}^NP(U_j=c_j)
=\prod_{j=1}^N\frac{1-\frac1{p_{j;A}}}{1-(\frac1{p_{j;A}})^k}(\frac1{p_{j;A}})^{c_j}\\
&=\frac1n\prod_{j=1}^N\big(1+\frac1{p_{j;A}}+(\frac1{p_{j;A}})^2+\cdots+(\frac1{p_{j;A}})^{k-1}\big)^{-1}.
\end{aligned}
\end{equation}

We have
\begin{equation*}
\log I_{N;A,2}=\sum_{j=1}^NU_j\log p_{j;A}.
\end{equation*}
Since
$$
\begin{aligned}
&EU_j=\frac{1-\frac1{p_{j;A}}}{1-(\frac1{p_{j;A}})^k}\sum_{m=0}^{k-1}m(\frac1{p_{j;A}})^m\\
&=\frac{1-\frac1{p_{j;A}}}{1-(\frac1{p_{j;A}})^k}\frac{1+(k-1)(\frac1{p_{j;A}})^k-k(\frac1{p_{j;A}})^{k-1}}{(1-\frac1{p_{j;A}})^2}\frac1{p_{j;A}},
\end{aligned}
$$
we have
\begin{equation}\label{EU}
EU_j\sim\frac1{p_{j;A}},\ \text{as}\ j\to\infty.
\end{equation}
From \eqref{ET} note that $EU_j$ and $ET_j$ have the same asymptotic behavior.
From \eqref{pjAasymp} and \eqref{EU}, we have
\begin{equation}\label{Elogagain}
E\log I_{N;A,2} \sim\sum_{j=1}^N\frac{\log p_{j;A}}{p_{j;A}}\sim\theta\log N,\ \text{as}\ N\to\infty.
\end{equation}
Note from \eqref{Elog} that $E\log I_{N;A,2}$ and $E\log I_{N;A,1}$ have the same asymptotic behavior.

We now give the appropriate redefinition of the mutually independent random variables
 $\{B_j\}_{j=1}^\infty$ and  $\{Xj\}_{j=1}^\infty$ that were defined in the proof of  part (i):

\noindent $B_j\sim \text{Ber}\big(\frac{\frac1{p_{j;A}}-(\frac1{p_{j;A}})^k}{1-(\frac1{p_{j;A}})^k}\big)$; that is,
\begin{equation}\label{qjagain}
q_j:=P(B_j=1)=1-P(B_j=0)=\frac{\frac1{p_{j;A}}-(\frac1{p_{j;A}})^k}{1-(\frac1{p_{j;A}})^k}.
\end{equation}

\noindent $X_j\stackrel{\text{dist}}{=}\log p_{j;A}\cdot U_j|\{U_j\ge1\}$;
that is
$$
P(X_j=m\log p_{j;A})=P(U_j=m|U_j\ge1)=\frac{1-\frac1{p_{j;A}}}{1-(\frac1p)^{k-1}}(\frac1{p_{j;A}})^{m-1},\ m=1,\cdots, k-1.
$$
As in part (i),
we have
$\mu_j:=EX_j\sim\log p_{j;A}$, as $j\to\infty$, and $\lim_{j\to\infty}\frac{X_j}{\mu_j}\stackrel{\text{dist}}{=}1$.
By the construction of $\{B_j\}_{j=1}^\infty$ and  $\{X_j\}_{j=1}^\infty$, we have
\begin{equation*}
\log I_{N;A,2}=\sum_{j=1}^NU_j\log p_{j;A}\stackrel{\text{dist}}{=}\sum_{j=1}^NB_jX_j.
\end{equation*}
By the same considerations as in part (i), it follows that
$W_N:=\frac{\sum_{j=1}^NB_jX_j}{E\sum_{j=1}^NB_jX_j}$ converges in distribution to $\frac1\theta D_\theta$; thus,
\begin{equation}\label{strangethm2}
\lim_{N\to\infty}\frac{\log I_{N;A,2}}{E\log I_{N;A,2}}\stackrel{\text{dist}}{=}\frac1\theta D_\theta.
\end{equation}

On the one hand, just as in \eqref{onehand}, by the convergence in distribution in \eqref{strangethm2} and the fact that the limiting distribution is a continuous one,   for any
 for any
sequence $\{\theta_N\}_{N=1}^\infty$ satisfying $\lim_{N\to\infty}\theta_N=\theta$, we have
\begin{equation}\label{onehandagain}
\lim_{N\to\infty}P(\frac{\log I_{N;A,2}}{E\log I_{N;A,2}}\le \frac1{\theta_N})=P(\frac1\theta D_\theta\le \frac1\theta)
=\frac{e^{-\gamma\theta}}{\Gamma(\theta+1)}.
\end{equation}
On the other hand,
let $\theta_N:=\frac{E\log I_{N;A,2}}{\log p_{N;A}}$ and note from \eqref{pjAasymp} and \eqref{Elogagain} that
$\lim_{N\to\infty}\theta_N=\theta$.
It follows from
 \eqref{1nprobagain} that
\medskip

$
P(\frac{\log I_{N;A,2}}{E\log I_{N;A,2}}\le \frac1{\theta_N})=P\big(I_{N;A,2}\le \exp(\frac{E\log I_{N;A,2}}{\theta_N})\big)=
P(I_{N;A,2}\le p_{N;A})
$
\medskip

$=\prod_{j=1}^N\big(1+\frac1{p_{j;A}}+(\frac1{p_{j;A}})^2+\cdots+(\frac1{p_{j;A}})^{k-1}\big)^{-1}
\sum^{'(k)}_{n\le p_{N;A}:p|n\Rightarrow p\in A  }\frac1n
$
\medskip

$=\frac{\sum^{'(k)}_{n\le p_{N;A}:p|n\Rightarrow p\in A}\frac1n}{\sum^{'(k)}_{n:p|n\Rightarrow p\in A_N}\frac1n}$.
From this and \eqref{onehandagain}, we conclude that
\begin{equation}\label{finalalmostagain}
\lim_{N\to\infty}\frac{\sum^{'(k)}_{n:p|n\Rightarrow p\in A_N}\frac1n}
{\sum^{'(k)}_{n\le p_{N;A}:p|n\Rightarrow p\in A}\frac1n}=e^{\gamma\theta}\Gamma(\theta+1),
\end{equation}
which is equivalent to part (ii) of Theorem \ref{GenMert}, just as \eqref{finalalmost} was equivalent to part (i) of the theorem. \hfill $\square$

\section{Proof of Theorem \ref{fixedN}}\label{3}
Fix $N\ge1$,  $k\ge2$ and $A_N=\{p_{1;A},\cdots, p_{N;A}\}$ as in the statement of the theorem.
Consider the random integer $I_{N;A,3}$ that was defined in \eqref{randomint3}, where $\{V_j\}_{j=1}^N$  are independent random variables with distributions given by
\eqref{Vdist}.
The support of  $I_{N;A,3}$ is the set of $k$-free integers in $\{n\in\mathbb{N}:p|n\Rightarrow p\in A_N\}$.
Let $f,g$ be functions satisfying
$$
\begin{aligned}
&f(j)=1, \ j=0,\cdots, k-2;\ \ f(k-1)=0;\\
&g(j)=j+1,\ j=0,\cdots, k-2;\ \  g(k-1)=k-1.
\end{aligned}
$$
For definiteness, we  take
$$
f(x)=1-\binom x{k-1},\ \ \ \  g(x)=(x+1)-\binom x{k-1}.
$$
Then  for arbitrary $n=\prod_{j=1}^Np_{j;A}^{c_j}$ in the support,
\begin{equation}\label{1key}
P(I_{N;A,3}=n)=\prod_{j=1}^NP(V_j=c_j)=\prod_{j=1}^N\frac{(p_{j:A}-1)^{f(c_j)}}{p_{j;A}^{\thinspace g(c_j)}}.
\end{equation}
Noting that $g(x)-f(x)\equiv x$, we rewrite the right hand side of \eqref{1key} as
\begin{equation}\label{2key}
\begin{aligned}
&\prod_{j=1}^N\frac{(p_{j:A}-1)^{f(c_j)}}{p_{j;A}^{\thinspace g(c_j)}}=\frac{\prod_{j=1}^N\big(\frac{p_{j:A}-1}{p_{j;A}}\big)^{f(c_j)}}{\prod_{j=1}^N p_{j;A}^{g(c_j)-f(c_j)}}\\
&=\frac{\prod_{j=1}^N\big(1-\frac1{p_{j;A}})^{f(c_j)}}{\prod_{j=1}^Np_{j;A}^{c_j}}=
\frac{\prod_{j=1}^N(1-\frac1{p_{j;A}})}{\prod_{j:c_j=k-1}(1-\frac1{p_{j;A}})}\times \frac1n\\
&=\Big(\prod_{j=1}^N(1-\frac1{p_{j;A}})\Big)\times\frac1{n_{\{(k-1)-\text{free}\}}}\times
\frac1{n_{\{(k-1)-\text{power}\}}\prod_{j:c_j=k-1}(1-\frac1{p_{j;A}})}\\
&=\Big(\prod_{j=1}^N(1-\frac1{p_{j;A}})\Big)\times\frac1{n_{\{(k-1)-\text{free}\}}}\times\frac1{n_{\{(k-1)-\text{power}\}}
\prod_{j:p_{j;A}|n_{\{(k-1)-\text{power}\}}}(1-\frac1{p_{j;A}})}\\
&=\Big(\prod_{j=1}^N(1-\frac1{p_{j;A}})\Big)\times\frac1{n_{\{(k-1)-\text{free}\}}}\times \frac1{\phi(n_{\{(k-1)-\text{power}\}})}.
\end{aligned}
\end{equation}
From \eqref{1key} and \eqref{2key} we obtain
\begin{equation}\label{keyfinal}
P(I_{N;A,3}=n)=\Big(\prod_{j=1}^N(1-\frac1{p_{j;A}})\Big)\thinspace \frac1{n_{\{(k-1)-\text{free}\}}\phi(n_{\{(k-1)-\text{power}\}})}.
\end{equation}
The theorem now follows from \eqref{keyfinal} along with the fact that
$$
\prod_{j=1}^N(1-\frac1{p_{j;A}})^{-1}=\sum_{n:p|n\Rightarrow p\in A_N}\frac1n
$$
and that
\medskip

$
\sum^{'(k)}_{n:p|n\Rightarrow p\in A_N}P(I_{N;A,3}=n)=1.
$
\hfill $\square$

\section{Proof of Theorem \ref{GenPhiMert}}\label{4}
Fix a subset $A\subset\mathbb{P}$ which satisfies $D_{\text{nat-prime}}(A)=\theta\in(0,1]$.
Denote the primes in $A$ in increasing order by $p_{1;A},p_{2;A},\cdots$, and let
$$
A_N=\{p_{1;A},p_{2;A},\cdots, p_{N;A}\}.
$$
Let $I_{N;A,3}$ be as in \eqref{randomint3},
where $\{V_j\}_{j=1}^\infty$  are independent random variables with distributions given by
\eqref{Vdist}.
We have
$$
\log I_{N;A,3}=\sum_{j=1}^N V_j\log p_{j;A}.
$$
It is easy to check that as with $ET_j$ and $EU_j$, we have
\begin{equation}\label{EV}
EV_j\sim\frac1{p_{j;A}},\ \text{as}\ j\to\infty.
\end{equation}
From
 \eqref{pjAasymp} and \eqref{EV}, we have
\begin{equation}\label{Elog3}
E\log I_{N;A,3} \sim\sum_{j=1}^N\frac{\log p_{j;A}}{p_{j;A}}\sim\theta\log N,\ \text{as}\ N\to\infty.
\end{equation}

As in the proof of Theorem \ref{GenMert}, we define mutually independent  random variables
 $\{B_j\}_{j=1}^\infty$ and  $\{Xj\}_{j=1}^\infty$:

\noindent $B_j\sim \text{Ber}(\frac1{p_{j;A}})$; that is,
\begin{equation}\label{qjagain}
q_j:=P(B_j=1)=1-P(B_j=0)=\frac1{p_{j;A}}.
\end{equation}

\noindent $X_j\stackrel{\text{dist}}{=}\log p_{j;A}\cdot V_j|\{V_j\ge1\}$;
that is
$$
P(X_j=m\log p_{j;A})=P(V_j=m|V_j\ge1)=\begin{cases}
(1-\frac1{p_{j;A}})(\frac1{p_{j;A}})^{m-1},\ m=1,2,\cdots, k-2.\\ (\frac1{p_{j:A}})^{k-2},\ m=k-1.\end{cases}
$$
As in the proof of Theorem \ref{GenMert}, we have
$\mu_j:=EX_j\sim\log p_{j;A}$ and $\lim_{j\to\infty}\frac{X_j}{\mu_j}\stackrel{\text{dist}}{=}1$.
By the construction of $\{B_j\}_{j=1}^\infty$ and  $\{X_j\}_{j=1}^\infty$, we have
\begin{equation*}
\log I_{N;A,3}=\sum_{j=1}^NV_j\log p_{j;A}\stackrel{\text{dist}}{=}\sum_{j=1}^NB_jX_j.
\end{equation*}
By the same considerations as in the proof of parts (i) and (ii) of Theorem \ref{GenMert}, it follows that
$W_N:=\frac{\sum_{j=1}^NB_jX_j}{E\sum_{j=1}^NB_jX_j}$ converges in distribution to $\frac1\theta D_\theta$; thus,
\begin{equation}\label{strangethm3}
\lim_{N\to\infty}\frac{\log I_{N;A,3}}{E\log I_{N;A,3}}\stackrel{\text{dist}}{=}\frac1\theta D_\theta.
\end{equation}

On the one hand, just as in \eqref{onehand} and \eqref{onehandagain}, by the convergence in distribution in \eqref{strangethm3} and the fact that the limiting distribution is a continuous one,   for any
sequence $\{\theta_N\}_{N=1}^\infty$ satisfying $\lim_{N\to\infty}\theta_N=\theta$, we have
\begin{equation}\label{onehandagainagain}
\lim_{N\to\infty}P(\frac{\log I_{N;A,3}}{E\log I_{N;A,3}}\le \frac1{\theta_N})=P(\frac1\theta D_\theta\le \frac1\theta)
=\frac{e^{-\gamma\theta}}{\Gamma(\theta+1)}.
\end{equation}
On the other hand,
let $\theta_N:=\frac{E\log I_{N;A,3}}{\log p_{N;A}}$ and note from \eqref{pjAasymp} and \eqref{Elog3} that
$\lim_{N\to\infty}\theta_N=\theta$.
It follows from
 \eqref{keyfinal} that
\medskip

$P(\frac{\log I_{N;A,3}}{E\log I_{N;A,3}}\le \frac1{\theta_N})=P\big(I_{N;A,3}\le \exp(\frac{E\log I_{N;A,3}}{\theta_N})\big)=
P(I_{N;A,3}\le p_{N;A})
$
\medskip

$=\prod_{j=1}^N(1-\frac1{p_{j;A}})\thinspace\sum^{'(k)}_{n\le p_{N;A}:p|n\Rightarrow p\in A}\thinspace \frac1{n_{\{(k-1)-\text{free}\}}\phi(n_{\{(k-1)-\text{power}\}})}$
\medskip

$=\frac{\sum^{'(k)}_{n\le p_{N;A}:p|n\Rightarrow p\in  A}\thinspace \frac1{n_{\{(k-1)-\text{free}\}}\phi(n_{\{(k-1)-\text{power}\}})}}
{\sum_{n:p|n\Rightarrow p\in A_N}\frac1n}$.
\medskip

\noindent From this and \eqref{onehandagainagain}  it follows that
$$
\lim_{N\to\infty}\frac{\sum_{n:p|n\Rightarrow p\in A_N}\frac1n}
{\sum^{'(k)}_{n\le p_{N;A}:p|n\Rightarrow p\in A}\thinspace \frac1{n_{\{(k-1)-\text{free}\}}\phi(n_{\{(k-1)-\text{power}\}})}}=e^{\gamma\theta}\Gamma(\theta+1),
$$
which is equivalent to \eqref{genphimert} just as
 \eqref{finalalmost} and \eqref{finalalmostagain} were equivalent to parts (i) and (ii) respectively  of Theorem \ref{GenMert} \hfill $\square$

\end{document}